\documentclass{article}
\usepackage[T2A]{fontenc}
\usepackage[cp1251]{inputenc}
\usepackage[english,russian]{babel}
\usepackage[tbtags]{amsmath}
\usepackage{amsfonts,amssymb,mathrsfs,amscd}
\usepackage{graphics}
\usepackage{euscript,textcomp,verbatim,fancyhdr}
\usepackage{latexsym}
\usepackage[hyper]{msb-a}
\JournalName{}
\theoremstyle{plain}
\newtheorem{Thm}{Theorem}
\newcommand{\RR}{\mathbb R}
\newcommand{\CC}{\mathbb C}
\newcommand{\NN}{\mathbb N}
\newcommand{\ZZ}{\mathbb Z}

\newcommand{\A}{\mathbf{A}}

\newcommand{\I}{\mathbf{I}}
\newcommand{\J}{\mathbf{J}}
\newcommand{\K}{\mathbf{K}}
\newcommand{\bS}{\mathbf{S}}
\newcommand{\T}{\mathbf{T}}
\newcommand{\m}{\mathbf{m}}
\newcommand{\z}{\mathbf{z}}

\newcommand{\mmu}{\boldsymbol{\mu}}
\newcommand{\Cr}{{\Sigma}}
\newcommand{\eps}{{\varepsilon}}

\newcommand{\cA}{\mathcal{A}}
\renewcommand{\d}{{\partial}}
\renewcommand{\hat}{\widehat}
\renewcommand{\tilde}{\widetilde}

\newcommand{\Int}{{\mathrm{int}\,}}

\newcommand{\rank}{\mathrm{rank\,}}
\renewcommand{\Re}{\mathrm{Re}}
\newcommand{\conv}{\mathrm{conv}}
\newcommand{\Q}{{C}}

\begin{document}

\title{The Liouville-type theorem for integrable Hamiltonian systems with incomplete flows}
\author{Elena A.~Kudryavtseva}
\address{Moscow State University}
\email{eakudr@mech.math.msu.su}
\date{}
\udk{}

\maketitle
\begin{fulltext}

\begin{abstract}
For integrable Hamiltonian systems with two degrees of freedom whose
Hamiltonian vector fields have incomplete flows, an analogue of the
Liouville theorem is established. A canonical Liouville fibration is
defined by means of an ``exact'' $2$-parameter family of flat
polygons equipped with certain pairing of sides. For the integrable 
Hamiltonian systems given by the vector field $v=(-\d f/\d w, \d f/\d z)$ on $\CC^2$
where $f=f(z,w)$ is a complex polynomial in $2$ variables, geometric properties of Liouville
fibrations are described.


\end{abstract}

\markright{Liouville theorem for integrable systems with incomplete flows}

\footnotetext[0]{This work was supported 
by the Government grant of the Russian Federation for support of research projects implemented by leading scientists, in the Federal State Budget Educational Institution of Higher Professional Education Lomonosov Moscow State University under the agreement No.~11.G34.31.0054.}

\section{Introduction} 

Suppose $(M^{2k},\omega)$ is a symplectic manifold and
$
F=(f_1,\dots,f_k):M^{2k}\to\RR^k
$ 
is a set of (functionally independent) smooth functions in involution. The triple
$(M^{2k},\omega,f_1)$ is called an {\it integrable Hamiltonian
system}, and the triple $L=(M^{2k},\omega,F)$ is called a {\it
Liouville fibration}. Two Liouville fibrations $L,L_1$ are called
{\it isomorphic} if $\omega=h^*\omega_1$ and $F=F_1\circ h$, for some
diffeomorphism $h:M^{2k}\to M_1^{2k}$. 

The classical Liouville theorem describes a Liouville fibration in a small neigh\-bo\-ur\-hood of
any compact connected regular fibre $T_\xi=F^{-1}(\xi)$ up to
isomorphism. The works \cite {Fl,BF,BC} stated the problem
of finding an analogue of the Liouville theorem for integrable
Hamiltonian systems whose Hamiltonian vector fields have incomplete
flows. We solve this problem for systems with $k=2$ degrees of
freedom (theorem \ref {thm:2}). We define a ``canonical'' Liouville
fibration $L_{W;U,d\I,\bS;[\Theta]}$ by means of so-called {\it
combinatorial-geometrical-topo\-logical data} (shortly {\it
CGT-data}), denoted by $(W;U,d\I,\bS;[\Theta])$, which are a
$2$-parameter family of flat polygons with certain pairing of sides.
Dynamical properties of the corresponding
Hamiltonian flow on an individual ``generic'' fibre $T_\xi$ were
studied by many authors (see \cite {konts,nov,buf} and references therein).

The author is grateful to A.\,T.~Fomenko for stating the problem,
L.~Bates, A.\,T.~Fomenko, T.\,A.~Lepskii, S.\,Yu.~Nemirovski,
A.\,I.~Sha\-farevich, I.\,Sypchenko, A.\,A.~Tuzhilin and A.\,B.~Zheglov
for useful discussions.

\section{Combinatorial-geometrical-topological data}

The CGT-data are defined as follows. 

\subsection{Combinatorial data}

Suppose $n\in\NN$, $\sigma\in\Sigma_{2n}$ is a permutation of the set $\{1,\dots,2n\}$, and 
$\tau:=\sigma^{-1}\in\Sigma_{2n}$ is the inverse permutation.
Denote the set of symbols $(a_1,a_1^{-1},\dots,a_n,a_n^{-1})$ by
$(w_{\tau(1)},\dots,w_{\tau(2n)})$. Thus
$w_{\tau(j)}=a_{k(j)}^{\eps(j)}$ where $k(j):=\lceil j/2\rceil$,
$\eps(j):=(-1)^{j-1}$, $1\le j\le2n$. A {\it quadratic word} is the
formal product
$$
 W:=w_1\dots w_{2n}=\prod_{j=1}^{2n}w_j
   =\prod_{j=1}^{2n}a_{k(\sigma(j))}^{\eps(\sigma(j))}.
$$

Consider the closed unit $2$-disk $D\subset\CC$ centred at the
origin, with $2n$ marked points $v_\ell=e^{\pi i\ell/n}\in\d D$,
$1\le\ell\le2n$. Consider the cell decomposition of $D$ into $2n$
vertices $v_1,\dots,v_{2n}$, $2n$ directed open edges and the open
$2$-cell $\Int D:=D\setminus\d D$. The edges are divided into pairs
$\alpha_k(t):=e^{\pi i(\tau(2k-1)+t)/n}$ and $\hat
\alpha_k(t):=e^{\pi i(\tau(2k)+1-t)/n}$ ($0\le t\le1$)
corresponding to the letters $a_k$ and $a_k^{-1}$ of $W$, $1\le k\le
n$. Denote by $\beta_k$ the directed curve on $D$ formed by the radii
that terminate at the centres of the edges $\alpha_k$ and
$\hat\alpha_k$, directed from $\alpha_k(1/2)$ to $\hat\alpha_k(1/2)$,
$1\le k\le n$. Consider the topological space
 $$
\Q=\Q_W:=D/\{\alpha_k(t)\sim\hat \alpha_k(t), \ 0\le t\le1\}_{k=1}^n
 $$
with quotient topology. It is a connected orientable closed surface
(of $\dim_\RR \Q=2$), with the induced cell decomposition having
$s:=n+1-2g$ vertices, $n$ directed edges identified with
$\alpha_1,\dots,\alpha_n$, and one $2$-cell. Here $g$ denotes the
genus of $\Q$. Denote by $\Q^r\subset \Q$ the $r$-skeleton of this
cell decomposition, $0\le r\le 2$. Thus
 $$
S=S_W:=\Q\setminus \Q^0
 $$
is a connected orientable genus-$g$ surface with $s$ punctures. If
the graph $\cup_{k=2g+1}^n\alpha_k$ is a spanning tree of the graph
$\Q^1$ (as can be achieved by renumbering the symbols
$a_1,\dots,a_n$) then the quadratic word $W$ will be called {\it
combinatorial data} of genus $(g,s)$.

\subsection{Geometrical data}

Suppose $U_1\subset\CC$ is an open subset, and
$$
\J=(J_1,\dots,J_n):U_1\to(\CC\setminus\{0\})^n, \qquad J_0:U_1\to\CC
$$
are continuous maps. The triple $(U_1,\J,J_0)$ is called {\it
geometrical data} with respect to the data $W$.

\subsection{Combinatorial-geometrical data}

Let us describe a natural geometrical object associated to the CG-data $(W;U_1,\J,J_0)$, namely a 2-parameter family of closed planar polygonal lines with certain pairing of sides.

The group $H_1(S)\cong\ZZ^n$ admits the basis
$\{[\beta_k]\}_{k=1}^n$. The group $H_1(\Q,\Q^0)\cong\ZZ^n$ admits
the bases $\{[\alpha_k]\}_{k=1}^n$ and
$\{[\tilde\alpha_k]\}_{k=1}^n$, where $\tilde\alpha_k:=\beta_k$ for
$1\le k\le2g$, $\tilde\alpha_k:=\alpha_k$ for $2g+1\le k\le n$. Here
$[\alpha_k], [\beta_k]$ denote the homology classes of the curves
$\alpha_k,\beta_k$ in appropriate homology groups. Consider the
coordinate isomorphisms 
 $
 \mmu=(\mu_1,\dots,\mu_n):H^1(S;\CC)\to\CC^n
 $ 
and
 $$
 \m=(m_1,\dots,m_n),\ \tilde\m=(\tilde m_1,\dots,\tilde m_n):H^1(\Q,\Q^0;\CC)\to\CC^n
 $$ 
with respect to these bases. Consider natural homomorphisms
 $$
 H_1(S) \stackrel{i}\longrightarrow
 H_1(\Q) \stackrel{p}\longrightarrow
 H_1(\Q,\Q^0).
 $$
Define the linear map $\A=(A_1,\dots,A_{2n}): \CC^n\to\CC^{2n}$ by
the rule
$$
A_\ell(\z):=\sum_{j=1}^{\ell-1}\eps(\sigma(j))z_{k(\sigma(j))}, \quad
 \z=(z_1,\dots,z_n)\in\CC^n,
 \ 1\le\ell\le2n.
$$
The homomorphism $i^*p^*:H^1(\Q,\Q^0;\CC)\to H^1(S;\CC)$ has
coordinate presentation $\T=(T_1,\dots,T_n)=\mmu
i^*p^*\m^{-1}:\CC^n\to\CC^n$ of the form
 $$
T_k(\m)
 =A_{\tau(2k)+1}(\m)-A_{\tau(2k-1)}(\m)
 =\sum_{\ell=1}^n\varphi_{\ell k}m_\ell,
 \quad 1\le k\le n.
 $$
Here $\varphi_{\ell k}:=\langle[\beta_\ell],[\beta_k]\rangle$ is the
intersection index of the cycles $[\beta_\ell]$ and $[\beta_k]$ in
$H_1(S)$. Thus $\rank(i^*p^*)=\rank(p^*)=2g$, and $\tilde
m_k=T_k(\m)$ for $1\le k\le2g$, $\tilde m_k=m_k$ for $2g+1\le k\le
n$.

Consider the $2$-parameter family of the closed $2n$-gonal lines on $\CC$ with consecutive vertices
$$
A_1(\J(\xi))+J_0(\xi),\dots,A_{2n}(\J(\xi))+J_0(\xi) \in\CC, \qquad \xi\in U_1,
$$ 
with parameter $\xi$, more precisely the $2$-parameter family of
the planar closed paths $\theta_{W;\J(\xi),J_0(\xi)}:\d D\to\CC$ formed by the segments
$$
 \theta_{W;\J(\xi),J_0(\xi)}(e^{\pi i(\ell+t)/n})
 =(1-t)A_\ell(\J(\xi))+tA_{\ell+1}(\J(\xi))+J_0(\xi), \quad 0\le t\le1, 
$$
$1\le\ell\le2n$, $\xi\in U_1$, where $A_{2n+1}:=A_1$.

\subsection{Topological data}

Take a point $\xi_0\in U_1$. A continuous map
$\Theta_{\xi_0}:D\to\CC$ is called {\it an extension} of the closed
$2n$-gonal line $\theta=\theta_{W;\J(\xi_0),J_0(\xi_0)}$ if
$\Theta_{\xi_0}|_{D\setminus\{v_\ell\}_{\ell=1}^{2n}}$ is an
orientation-preserving immersion and $\Theta_{\xi_0}|_{\d D}=\theta$.
So, an extension may have branching at the vertices $v_\ell$. In
general, an extension is not unique up to homeomorphisms of $D$
identical on $\d D$ (e.g.\ the octagon with consecutive vertices
$-2-i$, $i$, $2-i$, $-1/2$, $2+i$, $-i$, $-2+i$, $1/2$ has
inequivalent extensions \cite {Po}). Any extension $\Theta_{\xi_0}$
of $\theta$ can be included into a family of extensions
$\Theta_\xi:D\to\CC$ of the closed $2n$-gonal lines
$\theta_{W;\J(\xi),J_0(\xi)}$, $\xi\in U$, such that the map
$$
 \Theta=\Theta_{W;U,\J,J_0;\theta}:U\times D\to\CC, \qquad (\xi,d)\mapsto\Theta_\xi(d),
$$
is continuous, where $U\subset U_1$ is a small enough
simply-connected neighbourhod of $\xi_0$. Consider the set $[\Theta]$
of maps obtained by composing $\Theta$ with homeomorphisms of
$U\times D$ that are identical on $U\times\d D$ and preserve each
fibre of the projection $U\times D\to U$. The set $[\Theta]$ is
called {\it topological data} with respect to the data
$(W;U,\J,J_0)$.

\section{``Exact'' geometrical data and canonical Liouville fibrations}

Consider the topological 4-manifold
 $$
M_{W;U}:=U\times S_W=U\times S
 $$
and its open subsets 
 $$
 M_0:=U\times(\Q\setminus \Q^1)\approx U\times\Int D, \qquad 
 M_k:=U\times(p_{W}(U_k^+\cup U_k^-)), \quad 1\le k\le n.
 $$ 
Here $p_{W}:D\setminus\{v_\ell\}_{\ell=1}^{2n}\to S_{W}$ is the
projection, $U_k^+$ and $U_k^-$ are small disjoint neighbourhoods
of the open edges $\alpha_k$ and $\hat\alpha_k$ (respectively) in
$D\setminus\{v_\ell\}_{\ell=1}^{2n}$. Consider the following
immersions of these subsets into $U\times\CC\subset\CC^2$:
 $$
 (\xi,\Theta^0):M_0 \looparrowright U\times\CC, \qquad (\xi,\Theta^k):M_k\hookrightarrow U\times\CC, \quad 1\le k\le n,
 $$
where $\xi:M_{W;U}\to U$ is the projection,
$\Theta^0:=\Theta|_{M_0}$; the map $\Theta^k:M_k\to\CC$, $1\le k\le
n$, is defined by the rules 
$$
\Theta^k(\xi,p_{W}(d)):=\left\{
\begin{array}{ll} 
                 \Theta(\xi,d) & \mbox{if}\quad d\in U_k^+,\\
\Theta(\xi,d)-(0,T_k(\J(\xi))) & \mbox{if}\quad d\in U_k^-,
\end{array}\right.
$$
$(\xi,d)\in M_k$. Denote by $\cA$ the atlas on $M_{W;U}$ formed by all coordinate charts with local coordinates $(\xi,\Theta^k)$ on $M_k$, $0\le k\le n$.

\begin{Thm} \label {thm:1}
The following (``exactness'') conditions are equivalent: 

{\rm(i)} the atlas $\cA$ on $M_{W;U}$ is smooth, and the symplectic $2$-form $\{\Re(d\xi\wedge
d\Theta^k)\}_{k=0}^n$ on $(M_{W;U},\cA)$ is well-defined; 

{\rm(ii)} $((M_{W;U},\cA),\{\Re(d\xi\wedge d\Theta^k)\}_{k=0}^n,\xi)$ is a
Liouville fibration; 

{\rm(iii)} the 1-forms $\Re(T_k(\J(\xi))d\xi)$ with $1\le k\le2g$ (equivalently, with $1\le k\le n$) on $U$ are smooth and closed; 

{\rm(iv)} the surfaces $\{(\xi,T_k(\J(\xi)))\mid\xi\in U\}$ with $1\le k\le2g$
(equivalently, with $1\le k\le n$) are smooth and Lagrangian in
$(U\times\CC,\Re(d\xi\wedge d\Theta))$.
\end{Thm}

Suppose 
$$
\Re(T_k(\J(\xi))d\xi)=2\pi\, dI_k(\xi), \qquad \xi\in U, \quad 1\le k\le2g,
$$ 
for some smooth map $\I=(I_1,\dots,I_{2g}):U\to\RR^{2g}$.
Thus the condition (iii) of theorem \ref {thm:1} holds. Put
$$
 \bS:=(\K,J_0)\in C^0(U,\CC^s) \quad \mbox{where} \quad 
 \K:=(J_{2g+1},\dots,J_n)\in C^0(U,\CC^{s-1}),
$$ 
thus $(T_1\circ\J,\dots,T_{2g}\circ\J,\bS)=(\tilde\m\circ\m^{-1}\circ\J,J_0)$.
Let us associate the triple $(U,d\I,\bS)$ to such geometrical data
$(U,\J,J_0)$. By misuse of language, the triple $(U,d\I,\bS)$ will be
called {\it exact geometrical data with respect to $W$}.

By theorem~\ref {thm:1}, any {\it exact CGT-data}
$(W;U,d\I,\bS;[\Theta])$ determine a unique (up to equivalence)
Liouville fibration, which is denoted by $L_{W;U,d\I,\bS;[\Theta]}$
and called a {\it canonical Liouville fibration}.

\begin{Thm} \label {thm:2a}
For any pair of canonical Liouville fibrations $L_{W;U,d\I,\K,J_0;[\Theta]}$
and $L_{W;U,d\I^*,\K^*,J_0^*;[\Theta^*]}$ with the same data $(W;U)$,
the following condi\-tions are equivalent: 

{\rm(i)} for some $\xi_0\in U$, there exists an isomorphism $h:M_{W;U}\to M_{W;U}$ of these fibrations identical on the set $\{\xi_0\}\times(\Q_W^1\setminus \Q_W^0)$; 

{\rm(ii)} $d\I=d\I^*$ and $\K=\K^*$;
$\Theta_\xi=\Theta^*_\xi\circ h_\xi+J_0(\xi)-J_0^*(\xi)$ for some family of
homeomorphisms $h_\xi:D\to D$ identical on $\d D$, $\xi\in U$; the 1-form
$\Re((J_0(\xi)-J_0^*(\xi))d\xi)$ on $U$ is smooth and closed.
\end{Thm}

\section{The Liouville-type theorem}

For any Liouville fibration $(M^4,\omega,F)$, denote by $V_1,V_2$ the
Hamil\-tonian vector fields with the Hamiltonian functions $f_1,f_2$.
On each regular fibre $T_\xi:=F^{-1}(\xi)$, consider the flat
Riemannian metric $g_\xi$ inverse to the bi-vector field $(V_1\otimes
V_1+V_2\otimes V_2)|_{T_\xi}$. Denote by $\overline{T_\xi}$ the
completion of $T_\xi$ with respect to this Riemannian metric. 

\begin{Thm} \label {thm:2}
Suppose that a Liouville fibration $(M^4,\omega,F)$ is topo\-logically
locally trivial; moreover any its fibre $T_\xi$ is regular and
connected, the completion $\overline{T_\xi}$ of $T_\xi$ is compact and
$0<|\overline{T_\xi}\setminus T_\xi|<\infty$, $\xi\in
F(M^4)\subset\RR^2\cong\CC$. 

Then any point $\xi\in F(M^4)$ has a neighbourhood $U$ such that the Liouville fibration
$(F^{-1}(U),\omega,F)$ is isomorphic to a canonical Liouville
fibration $L_{W;U,d\I,\bS;[\Theta]}$, for some exact CGT-data of
genus $(g,s)$, where $s:=|\overline{T_\xi}\setminus T_\xi|$,
$g:=1-(\chi(T_\xi)+s)/2$. The flat Riemannian metric $g_\xi$ on any
fibre $T_\xi$ has a conical singularity at any puncture, where all
cone angles are integer multiples of $2\pi$.
\end{Thm}

\section{Examples via complex polynomials and their Newton polygons}

Consider a non-constant polynomial $f(z,w)=\sum_{\ell,m\ge0}a_{\ell
m}z^\ell w^m$ in two complex variables. H.~Flaschka \cite {Fl} and
A.\,I.~Shafarevich obser\-ved that the triple $(\CC^2,\Re(dz\wedge
dw),f)$ is a Liouville fibration. The set $\Cr_f$ of critical values
of $f$ is known to be finite \cite[\S2]{kl:metric}. The {\it
Newton polygon} of $f$ is the convex hull of its ``support'':
 $$
\Delta_f:=\conv\left\{(\ell,m)\in\ZZ^2\mid a_{\ell m}\ne0\right\} \subset \RR^2.
 $$ 
For any side $e$ of $\Delta_f$, denote
$f_e(z,w):=\sum_{(\ell,m)\in e}a_{\ell m}z^\ell w^m$. The polynomial
$f$ is called {\it weakly nondegenerate} with respect to $\Delta_f$
if, for any side $e$ of $\Delta_f$ that does not lie on coordinate
axes and for any point $(z,w)\in(\CC\setminus\{0\})^2$ with
$f_e(z,w)=0$, one has $df_e(z,w)\ne0$.

\begin{Thm} \label {thm:3}
Suppose a polynomial $f=f(z,w)$ is weakly nondege\-nerate with respect
to its Newton polygon $\Delta_{f}$, moreover
$\Delta_{1+z^\ell+w^m}\subseteq\Delta_f\subseteq\Delta_{(1+z^\ell)(1+w^m)}$
for some $\ell,m\in\NN$. Let $g>0$ be the number of integer points in
the interiour of $\Delta_{f}$, and $s-1$ be the number of integer
points of $\d\Delta_{f}$ that do not lie on coordinate axes.
Then the Liouville fibration $(\CC^2\setminus
f^{-1}(\Sigma_f),\Re(dz\wedge dw),f)$ satisfies all the hypothesis of
theorem {\rm\ref {thm:2}}. 

Moreover, the exact CGT-data $(W;U,d\I,\bS;[\Theta])$ of any corresponding canonical Liouville fibration have genus $(g,s)$ and satisfy the following: 

{\rm(i)} there exist functions $I_{2g+1},\dots,I_n,I_0\in C^\infty(U,\RR)$ such that
$\Re(J_0(\xi)d\xi)=dI_0(\xi)$, $\Re(J_k(\xi)d\xi)=dI_k(\xi)$,
$2g+1\le k\le n$, $\xi\in U$, where $\bS=(J_{2g+1},\dots,J_n,J_0)\in
C^0(U,\CC^s)$; 
moreover the functions $I_0,\dots,I_n\in C^\infty(U,\RR)$ are real parts of some holomorphic functions on $U\subset\CC$; 

{\rm(ii)} for any $\xi\in\CC\setminus\Cr_f$, there exists a bijection between the punctures of $T_\xi$ and the couples of neighbour integer points $A,B\in\d\Delta_f$ that do not
belong to the same coordinate axis, satisfying the following
condition: the cone angle of the flat Riemannian metric $g_\xi$ at
the puncture corresponding to $\{A,B\}$ equals $4\pi s_{\{A,B\}}$
where $s_{\{A,B\}}$ is the area of the triangle $A,B,(1,1)$.
\end{Thm}

Theorem \ref {thm:3} can be proved by using
\cite{kl:metric,hov.2,kl,kl:vta}. Its analogues for hyperelliptic
polynomials are given in \cite {L,kl,kl:vta}.

\end{fulltext}

\renewcommand {\refname}{References}


\begin{thebibliography}{99}

\bibitem {Fl} {\it Flaschka H.}
A remark on integrable Hamiltonian systems
// Physics Letters A. 1988. V.~131. \No~9. P.~505--508.

\bibitem {BF} {\it Bolsinov A.\,V., Fomenko A.\,T.}
Some actual unsolved problems on topology of integrable Hamiltonian systems
//
In book ``Topological methods in theory of Hamiltonian systems''. M.: Izd-vo
Faktorial, 1998. P.~5--23 (in Russian).

\bibitem {BC} {\it Bates L., Cushman R.}
Complete integrability beyond Liouville-Arnol'd
// Rep.\ Math.\ Phys. 2005. V.~56. \No~1. P.~77--91.

\bibitem{konts} {\it Kontsevich~M.\,L., Zorich~A.\,V.}
Connected components of the moduli spaces of Abelian differentials with prescribed singularities
// Inventiones mathematicae. 2003. V.~153. \No~3. P.~631--678.

\bibitem{nov} {\it Novikov~S.\,P.}
Topology of generic Hamiltonian foliations on Riemann surfaces
// Mosc. Math. J. 2005. Т.~5. \No~3. С.~633–-667.

\bibitem{buf} {\it Bufetov~A.\,I.}
Decay of correlations for the Rauzy-Veech-Zorich induction map on the space of interval exchange transformations and the central limit theorem for the Teichm\"uller flow on the moduli space of abelian differentials
// J.\ Amer.\ Math.\ Soc. 2006. V.~19. \No~3. P.~579--623.

\bibitem {Po} {\it Po\'enaru V.}
Extension des immersions en codimension one (d'apr\'es S.~Blank)
// S\'eminaire Bourbaki. 1968. \No~342. P.~473--505.

\bibitem{kl:metric}
{\it Kudryavtseva E.\,A., Lepskii T.\,A.}
Integrable Hamiltonian systems with incomplete flows and Newton's polygons
// Contemp.\ Probl.\ Mathem.\ Mechan. 2011. V.~VI. \No~3. P.~42--55 (in Russian).

\bibitem{hov.2}
{\it Khovanskii A.\,G.}
Newton polyhedra, and the genus of complete intersections
// 
// Funct.\ Anal.\ Appl. 1978. V.~12. \No~1. P.~38--46.

\bibitem{kl}
{\it Kudryavtseva E.\,A., Lepskii T.\,A.}
The topology of Lagrangian foliations of integrable systems with hyperelliptic Hamiltonian
// Sb.\ Math. 2011. V.~202. \No~3--4. P.~373--411.

\bibitem{kl:vta}
{\it Kudryavtseva E.\,A., Lepskii T.\,A.}
Topology of foliation and the Liouville theorem for integrable systems with incomplete flows
// Sb.\ Tr.\ Sem. Vekt. Tenz. Analizu. 2012. V.~XVII. P.~104--148 (in Russian).

\bibitem{L}
{\it Lepskii T.\,A.}
Incomplete integrable Hamiltonian systems with complex polynomial Hamiltonian of small degree
// Sb.\ Math. 2010. V.~201. \No~10. P.~1511--1538.

\end{thebibliography}
\end{document}